\documentclass[conference]{IEEEtran}
\IEEEoverridecommandlockouts
\usepackage[utf8]{inputenc}
\usepackage[T1]{fontenc}
\usepackage[english]{babel}
\usepackage{amssymb,amsfonts}
\usepackage{array}
\usepackage{mathtools}
\usepackage{amsmath}
\usepackage{enumerate}
\usepackage{enumitem}
\usepackage{comment}
\usepackage{caption}
\usepackage{subcaption}
\usepackage{hyperref}
\usepackage{bbm}
\usepackage{tikz}
\usetikzlibrary{plotmarks}
\usepackage{pgfplots}
\pgfplotsset{width=7cm,compat=1.8}
\usepackage{xcolor}
\usepackage{float}
\usepackage{verbatim}
\usepackage{eurosym}
\usepackage[noend]{algpseudocode}
\usepackage{filecontents}
\usepackage{algorithm}
\let\oldReturn\Return
\renewcommand{\Return}{\State\oldReturn}
\usetikzlibrary{plotmarks}
\usepackage[top=2cm, bottom=2cm, left=2cm, right=2cm]{geometry}
\DeclareMathOperator*{\argmin}{arg\,min}

\DeclareMathOperator{\prox}{prox}
\DeclareMathOperator{\card}{card}

\newtheorem{proposition}{Proposition}[section]

\def\K{\mathcal{K}}
\def\I{\mathcal{I}}

\newcommand\be{\begin{equation}}
\newcommand\ee{\end{equation}}

\makeatletter
\newcommand{\authorfootnotes}{\renewcommand\thefootnote{\@fnsymbol\c@footnote}}%
\makeatother

\makeatother

\renewcommand{\date}{\today}

\everymath{\displaystyle}

\begin{document}

\title{Computation and implementation of an optimal mean field control for smart charging}
\author{Adrien Séguret, Cheng Wan and Clémence Alasseur
\thanks{A. Seguret, C. Wan and C. Alasseur are with OSIRIS department in EDF Lab Saclay, France.}% <-this % stops a space
\thanks{A. Seguret is with CEREMADE, Université Paris Dauphine, France.}% <-this % stops a space
%\author{}
\thanks{This research benefited from the support of the FiME Lab (Institut Europlace de Finance) and the support of the FMJH Program Gaspard Monge for optimization and operations research and their interactions with data science}
}
\thispagestyle{plain}
\pagestyle{plain}

\maketitle

\begin{abstract}
This paper addresses an optimal control problem for a large population of identical plug-in electric vehicles (PEVs). The number of PEVs being large,  the mean field assumption is formulated to describe the evolution of the PEVs population and its interaction with the central planner. The resulting problem of optimal control of partial differential equations (PDEs) is discretized. Using convex analysis tools, we show the existence of an optimal solution and the convergence of the Chambolle-Pock algorithm to a solution. The implementation of this optimal control to the finite population of PEVs is detailed and we illustrate our approach with two numerical examples.
\end{abstract}

\IEEEpeerreviewmaketitle

\section{Introduction}

\subsection{Motivations}

The current energy transition that most countries are going through is  accompanied by a gradual but steady change in the means of transport. An increasing part of vehicles powered by fossil fuels are replaced by electric ones. This transformation raises new challenges as well as numerous benefits. On the one hand, the production and daily use of electrical vehicles (EVs) emit less pollution than fuel and diesel ones. The increasing capacity of batteries leads to considering fleet of EVs as non negligible virtual battery, the management of which is currently an important field of research. Plug-in electrical vehicles (PEVs) are expected to take part in ancillary services such as frequency and voltage regulation \cite{tomic2007using},\cite{zhong2014coordinated},\cite{wenzel2017real}, peak shaving \cite{uddin2018review}, valley filling \cite{gan2012optimal} and \cite{ma2011decentralized}, spinning reserve \cite{pavic2015value} and demand side management \cite{mets2011exploiting}.
On the other hand, it is well known that a non coordinated load of a large fleet of PEVs  can disrupt the energy balance on the electrical network during peak hours \cite{lopes2010integration}, \cite{deb2018impact}.
Several charging  scheme scenarios for a population of PEVs have been investigated over the last decade. 

\subsection{Related literature}

The decentralized and centralized architectures have been intensively inspected to optimize the scheduling of the charging of the fleet of PEVs. In a decentralized structure, it is common that either the PEVs interact, without any central planner, in a game setting \cite{couillet2012mean,ma2011decentralized}, or they optimize their consumption, in response to broadcasted signals (typically prices) by a central planner, to coordinate the consumption of the fleet \cite{gan2012optimal}.
In a centralized architecture, a central planner directly schedules and controls the charging of PEVs, with a view to taking part in ancillary services, as cited above. A wide range of tools from optimization has been explored to determine the optimal scheduling: linear and quadratic programming \cite{ayyadi2019optimal,mets2010optimizing}, dynamic programming \cite{lopez2018demand}, mixed-integer linear programming \cite{franco2015mixed}, stochastic optimization \cite{luo2017stochastic} etc. 
In spite of this diversity of techniques, centralized approaches may be computationally hardly scalable when the number of PEVs or the length of the planning time horizon are increasing. When the number of PEVs is large, one way to deal with the curse of dimensionality \cite{bellman2015adaptive} is to model the fleet of PEVs as a continuum of PEVs. this leads to a problem of optimal control of PDEs. It has been successfully applied in 
 \cite{bashash2011transport}, where the  aggregate charging power of PEVs is adapted to highly intermittent renewable power, and in \cite{ebrahimi2014aggregate} to perform renewable power tracking. The optimal control of PDEs for smart charging has been studied in a linear optimization framework in  \cite{le2015optimal},  where PEV fleets provide  regulation and vehicle-to-grid (V2G) services, and in \cite{sheppard2017optimal} where the optimal dispatch of an autonomous fleet of taxis is analyzed to serve
passengers and electric power demand during outages. The optimal control of PDEs is adopted in a linear quadratic framework  in \cite{le2016pde}, combined with Model Predictive Control (MPC) techniques, to solve the optimization problem of smart charging in real time. More recently, the authors solved numerically  a convex problem of optimal control of PDEs in \cite{seguret2021mean}, using the Chambolle-Pock Algorithm \cite{chambolle2011first}.

\subsection{Our contributions}

We provide a centralized optimal strategy to charge a large fleets of PEVs, over a finite time horizon.
While a large part of the litterature considers a continuous range of charging rates, we consider a finite number of charging rates ("not charging", "charging in fast or slow mode", etc.). It is because residential charging is mostly done at discrete rates \cite{nimalsiri2019survey}. The paper is based on the model and method developed in \cite{seguret2021mean}.
 
  The problem here is the optimal control of a PDE problem, which is discretized in time and space settings. 
We show the existence of a solution of the discretized problem, and give sufficient conditions for the convergence of the Chambolle-Pock Algorithm to the solution of the problem.
Among the existing literature, \cite{le2015optimal} and \cite{le2016pde} are the most close to our paper. However, our contributions are fairly different from those papers in several regards. Firstly we consider convex cost functions in contrast to linear and quadratic ones in  \cite{le2015optimal} and \cite{le2016pde}. Secondly constraints on the total mass per mode of charging is considered as well as constraints on the state of charging (SoC for short) of the PEVs at the final moment of the period, while \cite{le2015optimal} considers an aggregate constraint on the V2G (Vehicle-to-Grid) service in the linear case only. Whereas \cite{le2015optimal} and \cite{le2016pde}, we consider the mitigation of charging cycles per PEVs and the synchronization effects.
Indeed, the synchronization of PEVs can disrupt energy balance on the electrical network \cite{turitsyn2010robust}, while high charging frequency,
 is responsible for more intensive battery aging and degradation \cite{de2015impact, malhotra2016impact}. Finally, the optimal control computed on a basis of a continuum of PEVs, aims at being used to charge a finite number of PEVs. Thus, it is important to implement this control on a finite population of PEVs to check if the results are close to the one obtained in the mean field model. This work has not been done in \cite{seguret2021mean}, and we detail the implementation in this paper. This implementation is not trivial and we provide a condition on the mesh used for the implementation to obtain better results.

The remainder of this paper is organized as follows. Section \ref{cont_model} presents the evolution of the system in a continuous time and space setting. The problem is discretized in Section \ref{discrete_section} and we present in Section \ref{cp_section} the Chambolle-Pock algorithm to solve this problem. The implementation of the optimal control to a finite population of PEVs is detailed in Section \ref{App_alpha_method}. Section \ref{case_study_pres} is dedicated to two case studies and the analysis of the results.

% Bien citer nimalsiri2019survey qui justifie mes choix de chargement discret et d'éviter les multiple switch

\section{Aggregate control in the continuous framework}
\label{cont_model}

The purpose of this work is to determine a centraliezd control strategy for the charging of a large population of PEVs, over a finite period of time, in order to minimize a certain cost function while satisfying a set of constraints. The objective function does not need to be composed of the energy consumption cost only. It can also takes into considerations the aging of devices and the distance from a target level of the SoC at the end of the period. When the population of PEVs is large, combinatorial techniques as well as optimal control tools fail to solve the problem due to the curse of dimensionality \cite{bellman2015adaptive}. 
It is more suitable to work with the distribution of the states of the overall PEV population. Mean field control and PDEs control techniques can then be applied in such framework. A continuum of vehicles is considered over a period $[0,T]$. The population is assumed to be homogeneous, in the sense that the batteries of the PEVs share the same characteristic (battery capacity, charging rate, etc.). Besides, all the charging terminals are supposed to be identical with discrete charging rates. The state of each vehicle is composed of two variables: the SoC of its battery $s\in[0,1]$ and its charging mode $i\in \K $. The charging mode of a vehicle indicates its charging rate ("not charging", "charging in fast or slow mode", "V2G", etc.).
Since the charging terminals deliver discrete charging rates, $\K $ is a finite set and $I:=\card(\K )$, the cardinality of $\K$. For any $t\in[0,T]$, $i\in\K $ and $s\in[0,1]$, $m_i(t,s)$ stands for the density of vehicles at time $t$ with a SoC $s$ in mode $i$, and $b_i(s)$ for the instantaneous charging rate of their battery. If $b_i\geq0$ then we assume $b_i(1)=0$ and $b_i$ is non increasing and
if $b_i\leq0$ then we assume $b_i(0)=0$. We denote respectively by $\bar{D}_i$ and $\underline{D}_i$ the maximum and minimum total mass in mode $i$ allowed at time $k$.
Thus, for any $i\in\K$, $t\in[0,T]$, it holds:
\begin{equation}
\label{ineq cons}
\underline{D}_i(t)\leq
 \int_0^1m_i(t,s)ds\leq \bar{D}_i(t).
\end{equation}
The inequality \eqref{ineq cons} enables to take into account global constraints such as power consumption limit or minimum of available frequency response. 
In addition, we require that only a small proportion $\varepsilon>0$ of PEVs has a SoC inferior or equal to $\underline{s}\geq 0$ at the final moment $T$: 
\begin{equation}
\label{minsoc}
\sum_{i\in\K} \int_0^{\underline{s}} m_i(T,s) ds \leq \varepsilon. 
\end{equation}
 The transition intensity, $\alpha_{i,j}(t,s)$, is the instantaneous frequency with which vehicles with a SoC $s\in[0,1]$ at time $t\in[0,T)$ are transferred from the mode $i\in\K $ to another mode $j$. The function $\alpha_{i,j}$  is determined by the central planner, and is non negative:
 \begin{equation}
 \label{consalpha} 0\leq \alpha_{i,j}(t,s).
 %\leq \bar{\alpha}.
 \end{equation}
 By determining the transition intensity $\alpha$, the central planner controls how PEVs switch from one charging mode to another.
 The density $m_i$ evolves throughout time, due to the charge of the vehicles, the battery drain, and the transfers operated. No vehicles arrive or leave the considered population of PEVs during the period $[0,T]$. Using the conservation law, it holds:
\begin{equation}
\label{fk:0}
\begin{array}{l}
\partial_t m_i(t,s)+\partial_s(m_i(t,s)b_i(s))=  \\
-\sum_{j\in\K\setminus\{i\}}(\alpha_{i,j}(t,s)m_i(t,s)-\alpha_{j,i}(t,s)m_{j}(t,s)),
\end{array}
\end{equation}
with the initial condition for all $i\in\K$:
\begin{equation}
\label{inital_cond}
m_i(0,\cdot) = \bar{m}^0_i,
\end{equation}
where $\{\bar{m}^0_i\}_i$ is the initial distribution of the modes, satisfying: $\sum_{i\in\K} \int_0^1\bar{m}^0_i(s)ds =1$. 
The equation \eqref{fk:0} indicates how the distribution of SoC and charging mode of the PEV population evolves throughout time, governed by the controls $\{\alpha_{i,j}\}_{i,j}$ and the charging or discharging rate that they induce.
The reader is referred to  \cite{le2016pde} for a detailed deduction of \eqref{fk:0}.
For any $i,j\in\K $ with $i\neq j$, the function $E_{i,j}$ is defined by: $E_{i,j}:=\alpha_{i,j}m_i$.
This variable is introduced to make the equation \eqref{fk:0} linear w.r.t. $(m,E)$.
%function convex (defined below by equation \eqref{defJ}).
Thus, equation \eqref{fk:0} becomes:
\begin{equation}
\label{fk:1}
\begin{array}{l}
\partial_t m_i(t,s)+\partial_s(m_i(t,s)b_i(s))= \\
-\sum_{j\in\K\setminus\{i\}}E_{i,j}(t,s)-E_{j,i}(t,s).
\end{array}
\end{equation}
The term $E_{i,j}(t,s)$ in \eqref{fk:1} captures the instantaneous flow of vehicles from mode $i$ to mode $j$. The constraint \eqref{consalpha} becomes:
\begin{equation}
\label{borne_E}
0\leq E_{i,j}(t,s).
\end{equation}

The central planner's optimal charging strategy is obtained by the resolution of the following optimization problem:
\begin{equation}
\label{prob:primdual}
\begin{array}{l}
 \underset{m,E}{\inf}\,J(m,E)  \\
 \mbox{subject to }\,\eqref{ineq cons}, \eqref{minsoc}, \eqref{inital_cond},\eqref{fk:1} \mbox{ and }\eqref{borne_E},
\end{array}
 \end{equation}
where the objective function $J$ is defined by:
 \begin{equation}
 \label{defJ}
 \begin{array}{ll}
  J(m,E):= &  \int_0^T\sum_{i\in\K } \Omega_i(t,m_i(t,\cdot))dt\\
 & +\int_0^T\sum_{\underset{i\neq j}{i,j\in\K }}
\Theta_{i,j}(t,m_i(t,\cdot),E_{i,j}(t,\cdot))dt\\
  & +\sum_{i\in\K}\Gamma_i(m_i(T,\cdot)),
 \end{array}
\end{equation}
where $\Omega_i(t,\cdot)$, $\Theta_{i,j}(t,\cdot,\cdot)$ and $\Gamma_i$ are lower semi-continuous (l.s.c.) and convex functions for any $i,j\in\K$ and $t\in[0,T)$. The function $\Omega_i$ can typically represent the power consumption cost or the reward for frequency response, $\Theta_{i,j}$ the penalty for switching from mode $i$ to $j$ and $\Gamma_i$ the final cost penalizing the distance of the terminal state distribution from a given distribution. 
\section{Discrete model}\label{discrete_section}

The problem \eqref{prob:primdual} is discretized in space and time so that we can solve it numerically. We denote by $h>0$ the length of each step in SoC space, and by $\Delta t>0$ the length of each time step.
%Let $h>0$ and $\Delta t>0$ be 
They are taken in such a way that $N_h$ and $N_T$ are integers, with $N_h:=1/h$ and $N_T:=1/\Delta t$. For any $k\in\{0,\ldots,N_T\}$ and $\ell\in\{0,\ldots,N_h-1\}$,
 $\tilde{m}_i^{k,\ell}$ denotes the proportion of vehicles in mode $i\in\K$ with a SoC lying in $[\ell h, (\ell+1)h)$ at time $t_k:=k\Delta t$; $\tilde{\alpha}_{i,j}^{k,\ell}$
 denotes the transition intensity from mode $i$ to mode $j$ at time $t_k$; and $\tilde{E}_{i,j}^{k,\ell}$ denotes the associated flow, defined by: $\tilde{E}_{i,j}^{k,\ell}:=\tilde{m}_i^{k,\ell}\tilde{\alpha}_{i,j}^{k,\ell}$.
 
 The below equations \eqref{switch:alpha}-\eqref{neg_evo} approximate the PDE  dynamics \eqref{fk:1}. A splitting method is adopted: the reaction terms (i.e. transfers due to $\alpha$) are first taken into account, by introducing the variable $\tilde{m}^{k+\frac{1}{2},\ell}_i$ defined by:
\begin{equation}
\label{switch:alpha}
\tilde{m}^{k+\frac{1}{2},\ell}_i := \tilde{m}^{k,\ell}_i +\Delta t \sum_{j\in \K\setminus \{i\}} \tilde{E}_{j,i}^{k,\ell}-\tilde{E}_{i,j}^{k,\ell};
\end{equation}
then, the advection (the transport due to $\{b_i\}_{i}$) in each mode is estimated. An upwind scheme, commonly used in finite volume methods \cite{perthame2003equations}, is applied  differentiating the cases according to the sign, supposed constant, of $b_i$. Applying upwind scheme, if $b_i\geq0$, then:
\begin{equation}
\begin{array}{l}
\label{pos_evo}
 \frac{\tilde{m}^{k+1,\ell}_i-\tilde{m}^{k+\frac{1}{2},\ell}_i}{\Delta t}\\
+ 
\frac{b_i(x_{\ell+\frac{1}{2}})\tilde{m}^{k+\frac{1}{2},\ell}_i -
b_i(x_{\ell-\frac{1}{2}})\tilde{m}^{k+\frac{1}{2},\ell-1}_i} {h}\\
= 0,
\end{array}
\end{equation}

and if $b_i\leq0$:
\begin{equation}
\label{neg_evo}
\begin{array}{l}
 \frac{\tilde{m}^{k+1,\ell}_i-\tilde{m}^{k+\frac{1}{2},\ell}_i}{\Delta t}\\
+ 
\frac{
b_i(x_{\ell+\frac{1}{2}})\tilde{m}^{k+\frac{1}{2},\ell+1}_i - b_i(x_{\ell-\frac{1}{2}}) \tilde{m}^{k+\frac{1}{2},\ell}_i} {h}\\
= 0,
\end{array}
\end{equation}
 where $x_{\ell} :=(\ell+1/2)h$. 
 The initial condition is given by:
 \begin{equation}
 \label{init_cond}
 \tilde{m}_i^{0,\ell} = \bar{m}^{0,\ell}_i,
 \end{equation}
where  $\bar{m}^{0,\ell}_i:=\frac{1}{h}\int_{x_\ell -h/2}^{x_\ell + h/2}\bar{m}^0_i(s)ds$.

We define: $N_m := (N_T+1)N_hI$, $N_E := N_TN_hI(I-1)$, and $\tilde{N}:= N_hN_T I$.
 The constraints \eqref{switch:alpha}-\eqref{init_cond} being linear in $\tilde{m}$ and $\tilde{E}$, there exists a linear operators $\tilde{C}:\mathbb{R}^{N_m}\times \mathbb{R}^{N_E} \to \mathbb{R}^{N_m}$ such that the constraint 
 \begin{equation}
 \label{linear:cons}
 \tilde{C}(\tilde{m},\tilde{E}) = (\bar{m}^{0}, 0_{N_T\tilde{N}})
 \end{equation}
 is equivalent to \eqref{switch:alpha}-\eqref{init_cond}.

 To guarantee the non negativity of $\tilde{m}_i^{k+\frac{1}{2},\ell}$ (see Proposition \ref{prop_num_schemm})  in \eqref{switch:alpha}, the vectors $\{\tilde{E}_{i,j}\}_{i,j}$ are constrained by:
 \begin{equation}
 \label{cons_E}
 \sum_{j\in\K ,j\neq i}\sum_{\ell =0}^{N_h-1}\tilde{E}^{k,\ell}_{i,j}\Delta t-\tilde{m}^{k,\ell}_{i}\leq 0.
 \end{equation}
 The non negativity constraint on $E$ in continuous settings \eqref{borne_E} gives:
 \begin{equation}
 \label{borne_E_dis}
 0\leq \tilde{E}^{k,\ell}_{i,j}.
 \end{equation}
 The space and time steps $h $ and $\Delta t$ are taken to satisfy the  Courant–Friedrichs–Lewy condition  (see \cite{hirsch1990numerical}):
$\textstyle \underset{i\in\K }{\max} (\Vert b_i\Vert_\infty)\frac{\Delta t}{h}\leq 1$, which makes sure that $\tilde{m}_i^{k+1,\ell}$ is positive in \eqref{pos_evo} and \eqref{neg_evo}.

We have the following result for the numerical scheme, which ensures the non negativity of $\tilde{m}$ and mass conservation.

\begin{proposition}\label{prop_num_schemm}
If \eqref{cons_E} and \eqref{borne_E} are satisfied, then the mass is positive, i.e. $\tilde{m}\geq 0$, and it is conserved for any $k\in \{0,\ldots,N_T-1\}$, i.e.  $\textstyle\sum_{i\in \mathcal{K}} \sum_{\ell=0}^{N_h-1}\tilde{m}_i^{k,\ell}=\sum_{\ell=0}^{N_h-1}\bar{m}_i^{0,\ell}$.
\end{proposition}
\textit{Proof:}See Appendix \ref{section_appendix}.

 For any $i\in\K $ and $k\in\{0,\ldots,N_T\}$, the constraint \eqref{ineq cons} is discretized by:
\begin{equation}
\label{contraintes:etat}
\underline{D}_i^k\leq \sum_{j=0}^{N_h}\tilde{m}_i^{k,j}h\leq \bar{D}_i^k,
\end{equation}
Finally, the minimal SoC constraint at the final moment is expressed by:
\begin{equation}
\label{dis_smin}
\sum_{i\in\K}\sum_{\ell = 0}^{\underline{s}^h} \tilde{m}_i^{N_T,\ell} \leq \varepsilon,
\end{equation}
where $\underline{s}^h = \lceil\underline{s}/h\rceil $, $\lceil a\rceil$ being the largest integer smaller than or equal to a. We highlight that $\underline{s}$ is taken arbitrary small but strictly positive, to satisfy more easily the constraint qualification and the feasibility of the problem. However, we consider $\varepsilon = 0$ in practice. 
The objective function $J^h$ is defined by:
\begin{equation}
\label{def_Jh}
\begin{array}{ll}
 J^h(\tilde{m},\tilde{E}):=&\sum_{k=0}^{N_T-1}
 \sum_{i\in \K} \Omega^h_{i}(t_k, \tilde{m}^{k}_i)\Delta t \\
 & +\sum_{k=0}^{N_T-1}\sum_{\underset{i\neq j}{i,j\in\K }}\Theta_{i,j}^h(t_k,\tilde{m}^{k}_i,\tilde{E}^{k}_{i,j})\Delta t \\
&+\sum_{i\in\I} \Gamma_i^h(\tilde{m}^{N_T}).
\end{array}
\end{equation}
Functions $\Omega^h_i$, $\Theta^h_{i,j}$ and $\Gamma^h_i$ are respectively the discretization of $\Omega_i$, $\Theta_i$ and $\Gamma_i$.
For any $k\in\{0,\ldots, N_T\}$, $i,j\in I$, the function $\Theta^h_{i,j}(t_k,\cdot,\cdot)$ is supposed to be l.s.c. and convex w.r.t. both variables.
Similarly, for any $k\in\{0,\ldots, N_T\}$ and $i\in \mathcal{K}$, $\Omega_i^h(t_k,\cdot)$ and $\Gamma_i^h$ are supposed to be l.s.c. and convex.
The discretized optimization problem is:
\begin{equation}
\label{disc:pb}
\begin{array}{l}
 \underset{\tilde{m}\in\mathbb{R}^{N_m},\tilde{E}\in\mathbb{R}^{N_E}}{\inf}J^h(\tilde{m},\tilde{E}) \\
 \mbox{subject to }\, \eqref{linear:cons}-\eqref{dis_smin},
\end{array}
\end{equation}

\begin{proposition}
\label{existence_sol}
Problem \eqref{disc:pb} has a solution.
\end{proposition}

\textit{Proof: }From Proposition \ref{prop_num_schemm}, one has for any $i\in \mathcal{K},\,k\in \{0,\ldots, N_T\},\,\ell \in \{0,\ldots, N_h-1\}$ that $0\leq m_i^{k,\ell }\leq 1/h$. Then using inequalities \eqref{cons_E} and \eqref{borne_E_dis}, we deduce that $\tilde{E}$ is also bounded. Thus, any minimizing sequence of \eqref{disc:pb} has a converging subsequence, whose limit satisfies \eqref{linear:cons}-\eqref{dis_smin}, since the constraints are linear. Using that $J^h$ is l.s.c. we deduce that this limit is a solution of \eqref{disc:pb}.\hfill{$\square$}

\section{A primal-dual algorithm to solve the discretized optimal control problem}\label{cp_section}

The cost function $J^h$ is convex but not necessarily linear or quadratic.  This leads us to use the Chambolle-Pock algorithm, which shows good results on similar optimization problems in \cite{briceno2018proximal}. This algorithm  is a first-order primal-dual algorithm particularly suitable for non-smooth convex optimization problems with known saddle-point structure. In addition it enables to consider optimization problems with a large set of constraints, like the dynamic constraints described her.

\subsection{Chambolle-Pock algorithm}

In \cite{chambolle2011first}, a primal-dual algorithm is introduced to solve an optimization problem of the form:
\begin{equation}
\label{prob:prim}
\underset{y\in\mathbb{R}^N}{\min}\,\varphi (y) +\psi( y),
\end{equation}
and its associated dual:
\begin{equation}
\label{prob:dual}
\underset{\sigma\in\mathbb{R}^N}{\min}\,\varphi^\ast (-\sigma)
+\psi^\ast(\sigma),
\end{equation}
where
$\varphi:\mathbb{R}^N\to \mathbb{R}\cup \{+\infty \}$ and $\psi:\mathbb{R}^N\to \mathbb{R}\cup \{+\infty \}$ are proper, convex and l.s.c. functions.
Let $\hat{y}$ and  $\hat{\sigma}$ be respective solutions of \eqref{prob:prim} and \eqref{prob:dual}. 
We recall the definition of the proximal operator $\prox_f$ of a function $f$ at a point $x\in\mathbb{R}^N$:
\begin{equation}
\prox_{f}(x):=\underset{z\in\mathbb{R}^{N}}{\argmin}\,\frac{1}{2}\Vert z-x \Vert_2^2+f(z).
\end{equation}
One can easily check that:
\begin{equation}
\tilde{z}= \prox_{f}(x) \iff (x-\tilde{z})\in \partial f(\tilde{z}),
\end{equation}
where  $ \partial f(x)$ denotes the subgradient of $f$ at $x$, defined by: $\partial f(x):=\{r \in \mathbb{R}^N,\, f(z)\geq f(x)+ \langle r,z-x\rangle\}$, $\langle \cdot, \cdot \rangle$ being the usual scalar product in $\mathbb{R}^N$. If $(\hat{y},\hat{\sigma})$ is a saddle point, i.e. $\varphi (\hat{y}) +\psi( \hat{y})= \varphi^\ast (-\hat{\sigma})
+\psi^\ast(\hat{\sigma})$ then:
\begin{equation}
\left\{
\begin{array}{l}
-\hat{\sigma}\in\partial\varphi(\hat{y})\\
 \hat{y}\in\partial \psi^\ast(\hat{\sigma}),
\end{array}
\right.
\end{equation}
which is equivalent to, for any $\gamma>0$ and $\tau>0$:
\begin{equation}
\label{fixedpoint:primal}
\left\{
\begin{array}{l}
\prox_{\tau\varphi}(\hat{y}-\tau\hat{\sigma})=\hat{y}\\
\prox_{\gamma\psi^\ast}(\hat{\sigma}+\gamma\hat{y})=\hat{\sigma}.
\end{array}
\right.
\end{equation}
From \eqref{fixedpoint:primal}, problem \eqref{prob:prim} can be solved using the following primal-dual algorithm:
\begin{enumerate}
\item Fix $\theta\in[0,1]$, $\gamma>0$ and $\tau>0$ such that $\gamma \tau <1$.
\item Initialization: $(y^0,\tilde{y}^0,\sigma^0)\in\mathbb{R}^N\times \mathbb{R}^N\times \mathbb{R}^N$.
\item For any iteration $k$, knowing  $(y^k,\tilde{y}^k,\sigma^k)$, the following steps are executed:
    \begin{itemize}
    \item $\sigma^{k+1}=\prox_{\gamma\psi^\ast}(\sigma^k+\gamma \tilde{y}^k)$
    \item $ y^{k+1} = \prox_{\tau\varphi}(y^k-\tau\sigma^{k+1})$
    \item $ \tilde{y}^{k+1} =  y^{k+1} +\theta (y^{k+1}-y^k)$
    \end{itemize}
\end{enumerate}
%If $(\hat{y},\hat{\sigma})$ is saddle point, then 
Accordingly to \cite[Theorem 1]{chambolle2011first}, $(y^k,\sigma^k)$ converges to a saddle point $(y^\ast, \sigma^\ast)$ of \eqref{prob:prim} and \eqref{prob:dual}.

\subsection{Reformalulation of the discretized optimal control problem and convergence of the Chambolle-Pock algorithm}

The goal of this section is to formulate the constraint problem \eqref{disc:pb} to a problem of the form of \eqref{prob:prim}, and to show that the  Chambolle-Pock Algorithm converges to is solution. 
The functions $\varphi:\mathbb{R}^{N_m}\times \mathbb{R}^{N_E} \to \mathbb{R}$ and $\psi:\mathbb{R}^{N_m}\times \mathbb{R}^{N_E} \to \mathbb{R}$ are defined by:
\begin{equation*}
\label{def_varphi}
 \varphi(\tilde{m},\tilde{E}):=
 \left\{
 \begin{array}{ll}
 J(\tilde{m},\tilde{E}) & \mbox{if }(\tilde{m},\tilde{E})\mbox{ satisfy }\eqref{cons_E}-\eqref{dis_smin}\\
 +\infty & 
 \end{array}
 \right.
\end{equation*}
and
\begin{equation*}
\label{def_psi}
 \psi(\tilde{m},\tilde{E}):=
 \left\{
 \begin{array}{ll}
 0 & \mbox{if }(\tilde{m},\tilde{E})\mbox{ satisfy }\eqref{switch:alpha}-\eqref{init_cond}\\
 +\infty & 
 \end{array}
 \right.
\end{equation*}
The function $\psi$ constrains the variables $(\tilde{m},\tilde{E})$ to satisfy the discretized dynamics. The function $\varphi$ ensures charging requests and imposes system constraints.

The problem \eqref{disc:pb} is equivalent to:
\begin{equation}
\label{prob_ap:prim}
\underset{\tilde{m}\in\mathbb{R}^{N_m},\tilde{E}\in\mathbb{R}^{N_E}}{\min}\,\varphi (\tilde{m},\tilde{E}) +\psi(\tilde{m},\tilde{E}).
\end{equation}
%\begin{remark}
%Recall that the central planner is looking for an optimal control $\tilde{\alpha}_{i,j}^\ast$. While the Chambolle-Pock algorithm returns $(\tilde{m}^\ast,\tilde{E}^\ast)$, a solution of the finite dimensional problem \eqref{disc:pb}, the associated $\tilde{\alpha}^\ast_{i,j}$ can be computed using the definition of $\tilde{E}$.
%\end{remark}
One needs to show that there exists a saddle point associated to the convex problem \eqref{prob_ap:prim} to ensure the convergence of the Chambolle-Pock Algorithm to the solution of \eqref{prob_ap:prim}. It is enough to show that strong duality condition is satisfied. From Proposition \ref{existence_sol}, we know that the infinimum of \eqref{prob_ap:prim} is attained. In practice, strong duality condition is satisfied if the period $T$ is not too short and if the lower bound  $\bar{D}_i$, for at least least one mode $i\in \mathcal{K}$ whose charging rate $b_i$ is non negative, is not too low. 
Nevertheless, we state in Appendix \ref{slater_suff_cond} a proposition giving a sufficient condition to satisfy the Slater’s constraint qualification \cite{boyd2004convex}.

Finally, next proposition states the convergence of the Chambolle-Pock algorithm when applied to problem \eqref{prob_ap:prim}.

\begin{proposition}
If problem \eqref{prob_ap:prim} has a finite solution and if strong duality property holds, then there exists a solution $\hat{\sigma}$ of the dual of problem of \eqref{prob_ap:prim}. In addition, the sequence $(y^k,\sigma^k)$, defined by the Chambolle-Pock Algorithm converges to a saddle point 
%$(\tilde{m}^\ast,\tilde{E}^\as)$ 
of \eqref{prob:prim} and \eqref{prob:dual}
\end{proposition}
\textit{Proof: }Existence of a saddle point is a direct consequence of the existence of a finite solution of \eqref{prob_ap:prim} and of strong duality \cite{boyd2004convex}. Since the convex problem  \eqref{prob_ap:prim} is finite dimensional and has a saddle point and, the convergence of  $(y^k,\sigma^k)$ to a saddle point is guaranteed by \cite[Theorem 1]{chambolle2011first}.\hfill{$\square$}

\section{Implementation of the optimal control to a finite population of vehicles}
\label{App_alpha_method}

This section aims at proposing a method to implement a mean field control $\alpha$ 
%(defined in Section \ref{cont_model})
to a population of $n$ PEVs. We assume that $\alpha$ is a continuous interpolation of the discrete optimal control $\tilde{\alpha}$ obtained by the Chambolle-Pock Algorithm.
We consider another time and space discretization parameters than the ones defined in \eqref{def_Jh}, to implement the control. We denote by $u>0$ the space step, $N_u$ the associated number of space intervals to cover $[0,1]$, and $\Delta r>0$ the time step in this section. The definition of another mesh is motivated by the following condition on $u$ and $\Delta r$:
\begin{equation}
\label{ineq_n_mesh}
\frac{T}{nu\Delta r}+\Vert \alpha \Vert_\infty( \Delta r+u) \ll 1.
\end{equation}
If the control $\alpha$ is continuous and implemented as it is proposed in this section, then the inequality \eqref{ineq_n_mesh} ensures that the empirical distribution of the vehicles is close to the solution of \eqref{fk:0} in the weak sense (this result is to be shown in a coming paper by one of the author), that we omit to specify here. 
This condition allows to ensure that the mean field assumption of a continuum of PEVs is relevant. 
If the optimal control $\tilde{\alpha}$ is computed on a fine mesh, then $n$ may not be large enough to satisfy \eqref{ineq_n_mesh}, if the same mesh is used for the implementation. To this end, it may be interesting to implement the optimal control on a coarser mesh.

%We recall that the functions $\alpha$ is positive and defined for any $t\in[0,T)$, $s\in[0,1]$ and $i,j\in \mathcal{K}$.
For any $z\in\{1,\ldots,n\}$ and  $t\in[0,T]$, the $z^{th}$ PEV is described by its state variable $x^z(t):=(s^z(t), q^z(t))\in[0,1]\times \mathcal{K}$, where $s^z(t)$ is the SoC and $q^z(t)$ the charging mode. For any $z$, the dynamic of $s^z(t)$ is defined by:
\begin{equation}
\label{evo_soc}
ds^z(t) = b_{q^z(t)}(s^z(t))dt,
\end{equation}
where $b_{q^z(t)}$ stands for the charging rate determined by the mode $q^z(t)\in\K$ before the vehicle's mode is changed to another one. 

At any time step $t_k\in\{\Delta r,\ldots, T\}$, the integer $a^n_{i,j}(k,\ell)$ corresponds to the number of vehicles with a SoC around $x_\ell$ and in mode $i$ to be transferred to mode $j$. It is defined by:
\begin{equation}
\label{determination_a}
a^n_{i,j}(k,\ell) : = \left\{
\begin{array}{ll}
\lfloor 
N_i^{k,\ell} A_i^{k,\ell}\Delta r
\rfloor
 & \mbox{if }\bar{D}_j<1, \vspace{0.2cm} \\
 \lceil 
N_i^{k,\ell} A_i^{k,\ell} \Delta r
\rceil
 & \mbox{otherwise},
\end{array}
\right.
\end{equation} 
where
\begin{equation}
A_i^{k,\ell}: = \frac{1}{u}\int_{x_{\ell}-u/2}^{x_{\ell}+u/2}
\alpha_{i,j}(t_k,s)ds,
\end{equation}
and $N_i^{k,\ell}$ is the number of vehicles at time $t_k$ in mode $i$ with a SoC lying in the range $[y_\ell-u/2,y_\ell+u/2)$, and  $ \lfloor a \rfloor $ is the largest integer smaller than or equal to a.
 At each time step $t_k\in\{0,\ldots,T-\Delta r \}$, the following steps are applied:\newline
\textbf{STEP 1}: For vehicles of any mode $i\in \K$ and SoC index $\ell\in \{0,\ldots,N_u-1\}$, transfers to the other modes are performed:
\begin{itemize}
\item Evaluation of $V_i^{k,\ell}$, the set of PEVs in mode $i$ with a SoC in the range $[y_\ell-u/2,y_\ell+u/2)$ at time $t_k$. Then, we define $N_i^{k,\ell}:= \card(V_i^{k,\ell})$. 
\item For any $j\in\K\setminus\{i\}$, evaluation of $a^n_{i,j}(k,\ell)$ using \eqref{determination_a}. The transferred PEVs are chosen on a first-in-first-out basis. We denote by $T_{i,j}^{k,\ell}$ the indices of the vehicles transferred. It holds: $T_{i,j}^{k,\ell}\subset V_i^{k,\ell}$. \item Update of the modes of all the transferred vehicles: $\forall z \in T_{i,j}^{k,\ell}$, $q^z(t_{k+1}): = j$. If the upper or the lower limit of number of vehicles for a certain mode $j$ is attained (see inequality \eqref{ineq cons}), transfers to this mode are stopped if the upper bound is attained, and transfers from this class are stopped if the the lower bound is attained.
\end{itemize}
\textbf{STEP 2}: the SoC of all the vehicles are updated according to \eqref{evo_soc}.

\section{Case study}\label{case_study_pres}
This section illustrates the model and furnishes numerical results of the implementation of an estimate of the optimal control, obtained after sufficient iterations of the Chambolle-Pock Algorithm. 
Two modes of charging are considered $\K:=\{0,1\}$. A vehicle receives zero charging power when being in mode $0$, and $P_{ON}>0$ when being in mode $1$. The central planner operates on a large population of PEVs. She decides at different time steps the proportion of PEVs to charge at zero charging power or at $P_{ON}$, depending on their SoC. There are no arrivals and departures of vehicles during this period. We assume that the central planner knows the SoC of each vehicle at time $0$.
Two different examples are exposed in this section. We highlight that the two examples represent two very different situations of smart charging, and that the corresponding optimization problems also have different structures.

The following parameters are chosen for the simulations.
The considered period is $T=5$h, the time and space steps are $\Delta t = 450$ seconds and $h = 0.05$. The charging rate is: $b_1(s)=1/45000$ for all $ s\in [0,0.75)$ and $b_1(s)=(1-s)/11250$ for all $s\in [0.75,1]$.
To take into account the battery drain for non charging vehicles, we have set $b_0(0)=0$ and $b_0(s)= 3.86\times 10^{-7} $ for any $s\in (0,1]$.
The delivered charging power is $P_{ON}=20$kW and is independent of the SoC $s$.

\subsection{The problem of response to peak and off peak pricing}\label{case 0}
The central planner operates on a population of $n=500$ PEVs located on the same node of the distribution grid. The site has a non flexible and deterministic electrical demand due to other (non-EV) electric usage: it may include some PEV charging that are not controlled by the central planner. A fixed controlled population of PEVs is plugged in this area over a period of $[0,T]$. The central planner takes part in ancillaries services by contracting a contract where the price of electricity depends on whether it is during  peak or off-peak hours. 
It is well known that peak and off-peak hours pricing  structure  can  potentially  shift  the peak load demand and thus perform peak-shaving and valley-filling  \cite{dubey2015determining}, \cite{dubey2015electric}. The increase of the tariff during the peak consumption period and the reduction of the tariff during low consumption period
induce users to charge during low load period so as to realize peak shift. The objective of this example is to understand how central planner responses to the price while satisfying several conditions: a minimal SoC at the end of the period and the limited transformer capacities.
The  synchronization of PEV is another phenomenon to avoid \cite{turitsyn2010robust}, which can disrupt energy balance on the electrical network. Another risk is the aging of the battery caused by  switching frequently between different modes of charging. Indeed the charge/discharge cycles directly damage the lifetime of a battery \cite{abronzini2019cost}. To deal with both problems, the functions $\Theta_{1,0}$ and $\Theta_{0,1}$ are adopted to penalize the proportion of PEVs whose charging mode are switched from $0$ to $1$ and from $1$ to $0$ respectively, with:
\begin{equation}
 \label{def_theta01}
\Theta_{i,j}(a,b):=
\left\{
\begin{array}{ll}
 \theta_{i,j}b+ \tilde{\theta}_{i,j}\frac{b^2}{a} & \mbox{if }a>0\mbox{ and }b\geq 0, \\
 0 & \mbox{if }a=0 \mbox{ and }b=0, \\
 +\infty&\mbox{otherwise},
\end{array}\right.
\end{equation}
where $\theta_{i,j},\tilde{\theta}_{i,j}\geq 0$ and $(i,j)\in\{(1,0),(0,1)\}$. Such functions are known to be convex and l.s.c. \cite{benamou2000computational}.
On the one hand, 
the second order term of $b$ in \eqref{def_theta01} penalizes large values of $\tilde{E}$  and thus large displacement from one mode to another. On the other hand, the linear part of $b$ in \eqref{def_theta01} penalizes small values of $\tilde{E}$ close to 0, which prevents multiple small switches between the modes, when the time horizon is long. When multiple small changes occur, PEVs tend to have many charge/discharge cycles. The function $\Omega^h_{1}$ corresponds to the running cost due to electricity consumption. The quantity $\textstyle \tilde{p}^k: = \int_{t_k}^{t_{k+1}}p(t)dt$ is the average cost of electricity (in \euro/kWh) over the period $[t_k,t_{k+1}]$. We define $\Omega^h_{1}$ by:
\begin{equation*}
\Omega^h_{1}(t_k,\tilde{m}_1^{k,\ell}):=P_{ON}\tilde{p}^k\tilde{m}^{k,\ell}_1.
\end{equation*}
Only charging vehicles consume electricity. 
Finally, the central planner minimizes the following cost:
\begin{equation}
\label{def_Jh_case_0}
\begin{array}{l}
 J^h(\tilde{m},\tilde{E}):=\\
\sum_{k=0}^{N_T-1}\sum_{\ell=0}^{N_h-1}
  \big(  \Theta_{0,1}(\tilde{m}^{k,\ell}_0,E_{0,1}^{k,\ell})
  +\Theta_{1,0}(\tilde{m}^{k,\ell}_1,E_{1,0}^{k,\ell})\big)
  h\Delta t\\
   +\sum_{k=0}^{N_T-1}\sum_{\ell=0}^{N_h-1} \Omega_{1}(t_k,\tilde{m}_1^{k,\ell}) h\Delta t
\end{array}
\end{equation}
The system must satisfy the following constraints:
\begin{itemize}
\item Constraint on the aggregate power: all the PEVs are connected at the same node in the distribution grid. As a consequence, the instantaneous total power delivered to the population of PEVs is bounded over the time. This power limit constraints the proportion of PEVs in mode $1$. For any $k\in\{0,\ldots,N_T \}$, the central planner must satisfy:
\begin{equation}
\label{contraintes:etat0}
\sum_{j=0}^{N_h}\tilde{m}_1^{k,j}h\leq \bar{D}_1,
\end{equation}
where $\bar{D}_1$ is the maximum proportion of vehicle that can be in mode $1$ over the period $[0,T]$ to respect the constraint of maximum power consumption.

\item Constraint on the final SoC of the PEVs: contrary to the model exposed in Section \ref{discrete_section}, we take $\varepsilon=0$, meaning that all the SoC must have attained a minimum $\underline{s}$ at the end of the period. Thus the central planner must satisfy:
\begin{equation}
\label{dis_smin0}
\sum_{\ell = 0}^{\underline{s}^h} \tilde{m}_0^{N_T,\ell} + \tilde{m}_1^{N_T,\ell} =0,
\end{equation}
\end{itemize}
It is assumed that the central planner knows in advance:
\begin{itemize}
\item The price of electricity $p(t)$ during the period $[0,T]$.
\item The upper limit of consumption $\bar{D}_1^k$ for any time interval $[t_k,t_{k+1})$.
\end{itemize}
The optimization problem is:
\begin{equation}
\label{disc:pb:0}
\begin{array}{l}
 \underset{\tilde{m}\in\mathbb{R}^{N_m},\tilde{E}\in\mathbb{R}^{N_E}}{\inf}J^h(\tilde{m},\tilde{E}) \\
 \mbox{subject to }\, \eqref{linear:cons}-\eqref{dis_smin},\eqref{contraintes:etat0}\mbox{ and }\eqref{dis_smin0},
\end{array}
\end{equation}
where $J^h$ is defined by equation \eqref{def_Jh_case_0}. 

The price of electricity $p$ for the simulation is displayed in Figure \ref{fr_comp}. The entire time period is composed of four different periods: two periods of peak hours and two periods of off-peak hours. The parameters for the switching costs are: $\theta_{1,0}=\theta_{0,1} = 0.04$ and $\tilde{\theta}_{1,0}=\tilde{\theta}_{0,1} =  20$. The minimum SoC at the end of the period is $\underline{s}=0.7$ and  the maximum proportion of vehicles that can be in mode 1 is constant over the time: $\bar{D}_1 = 1/3$.
The Chambolle-Pock algorithm has performed 12000 dual primal iterations with parameters $\theta = \gamma = 0.5$ and $\tau =1.8$.
 An approximation of the optimal discrete variables $(m,E)$ of the problem \eqref{disc:pb:0} is obtained and its associated $\alpha$ is then determined using $\alpha_{i,j}^{k,l} = E_{i,j}^{k,l}/m_i^{k,l}$. 
The control $\alpha$ is applied to the finite population of PEVs, as described in Section \ref{App_alpha_method}, using the same space and time mesh as for the Chambolle-Pock algorithm. The number of PEVs per charging mode over the period $[0,T]$ is displayed in Figure \ref{empirical_nb_vec}. By comparing with the curve of price in Figure  \ref{fr_comp}, one observes that the vehicles  charge only during the off-peak hours, when the price is lower. The first off-peak period starts at time $t = 1$h: at that moment the proportion of vehicles in mode 1 starts to increase progressively and
reaches its limit $\bar{D}_1$, before the mid-term of the first off-peak period. The decrease of the number of PEV in mode 1 starts before the end of the first off-peak period and is progressive as well. The second off-peak period, starting at time $t=3.5h$,  is handled exactly in the same way. All the PEVs start and finishe in mode 0, except 6 of them which finish in mode 1.
The total number of changes from one mode to another among the PEV population during the whole period is 927. More precisely, 33 PEVs have not changed at all their mode of charging, 6 PEVs have changed once  and 461 PEVs have changed twice. 
Thus, the PEVs that are charged during the first off-peak period are different from those being charged during the second off-peak period. Figure \ref{first_time_soc_07} represents the distribution of the time, when the SoC attains $0.7$. It shows that a large proportion of PEVs SoC reaches $0.7$ during the first off-peak period, between $t=1$ and $t=2$. One can also observe that, after $t=4$, no PEV reaches SoC 0.7 during the last hour of the period.
The penalization on the transfers helps to avoid both excessive changes between different charging modes per PEV and the synchronization of PEVs. Indeed, as observed in the Figure \ref{empirical_nb_vec}, the increase and decrease of the number of charging PEVs are always progressive. In addition, all the PEVs do not stop charging when their SoC reaches $0.7$. As observed in Figure \ref{final_initial_distrib}, most of PEVs have a SoC greater than 0.7 at the end. 

The Figure \ref{theorical_nb_vec} represents the proportion of PEVs in each mode over the time, when the control is applied to a continuum of PEV.
Figures \ref{theorical_nb_vec} and \ref{empirical_nb_vec} are very similar, where Figure \ref{theorical_nb_vec} represents the result with a continuum of PEVs.
Observing that Figures \ref{theorical_nb_vec} and \ref{empirical_nb_vec} are very similar, this highlights that the procedure, described in Section \ref{App_alpha_method}, to implement the control $\alpha$ on a finite number of PEVs, is able to obtain an empirical distribution of PEVs close to the controlled distribution when a continuum is considered. 
The initial and final empirical distribution of the SoC of the fleet of PEVs are displayed in Figure \ref{final_initial_distrib}. At time $t=0$, the SoC are bounded between 0.2 and 0.6 while at time $t=5$h, almost all the SoC are at level 0.7 as desired. One can observe few vehicles whose SoC is lower than 0.7. The final constraint is not satisfied for 4.6\%  of PEVs here. This effect is expected: the control  $\alpha$ applied to the finite population was computed in a mean field framework. A simple heuristic can be used in a post processing step to take over vehicles whose final SoC are lower than 0.7. The more important number of PEVs is, the better results the mean field control $\alpha$ shows. For instance, a proportion of 10\% of the PEVs does not satisfy the final constraint for the SoC, when there are only 100 PEVs. This proportion is divided by two in our simulation with 500 PEVs.
In Figure \ref{evolution_soc}, the evolution of the SoC of thirty PEVs is displayed over the period $[0,T]$ . As expected, the SoC at the final time are equal or superior to 0.7. Also, one can see that the PEVs charge only during the off-peak period, i.e. between $t = 1$h and $t = 2$h, and after $t = 3$h30. 

\begin{figure*}[t]
\centering
\begin{minipage}[t]{0.32\linewidth}
\centering
\begin{tikzpicture}[scale=0.75, transform shape]
\begin{axis}[
xlabel={Time}, ylabel = {Proportion of vehicles}, legend entries={Mode 0, Mode 1}, ymin=0, ymax =1, xmin=0,   xmax=5,
legend style={at={(1,0)},anchor=south east}, ytick={0,0.10, 0.2,0.3,0.4,0.5,0.6,0.7,0.8,0.9,1},
 bar width=2pt,
 ymajorgrids=true,
ybar stacked
]
\addplot table[ybar interval=0,x index=0,y index=1]{mass_mode_theoritical_case_1.txt};
\addplot table[ybar interval=0,x index=0,y index=2]{mass_mode_theoritical_case_1.txt};
\end{axis}
\end{tikzpicture}
\captionsetup{font=scriptsize}
\caption{Proportion of vehicles per mode of charging over the period when a \textit{continuum} of vehicles is considered}
\label{theorical_nb_vec}.\end{minipage}
\begin{minipage}[t]{0.32\linewidth}
\begin{tikzpicture}[scale=0.75, transform shape]
\begin{axis}[
xlabel={Time}, ylabel = {Proportion of vehicles}, legend entries={Mode 0, Mode 1}, ymin=0, ymax =1, xmin=0,   xmax=5,
legend style={at={(1,0)},anchor=south east}, ytick={0,0.10, 0.2,0.3,0.4,0.5,0.6,0.7,0.8,0.9,1},
 bar width=2pt,
 ymajorgrids=true,
ybar stacked
]
\addplot table[ybar interval=0,x index=0,y index=1]{mass_mode_case_1.txt};
\addplot table[ybar interval=0,x index=0,y index=2]{mass_mode_case_1.txt};
\end{axis}
\end{tikzpicture}
\captionsetup{font=scriptsize}
\caption{Proportion of vehicles per mode of charging over the period when the control is implemented to a \textit{population of 500PEVs}.}
\label{empirical_nb_vec}
\end{minipage}
\begin{minipage}[t]{0.32\linewidth}
\centering
\begin{tikzpicture}[scale=0.75, transform shape]
\begin{axis}
[ymax = 0.18,
xmin = 0, xmax = 5,legend style={at={(1,1)},anchor=north east}, 
xtick  = {0,1,2,3,4,5},
grid=major,
xlabel={Time (h)},
ylabel={Price \euro/kWh}]
\addplot [draw=blue] table[x index=0,y index=1]{Price_case_1.txt};
\end{axis}
\end{tikzpicture}
\captionsetup{font=scriptsize}
\caption{Price of electricity over the period}
\label{fr_comp}
\end{minipage}

\begin{minipage}[t]{0.32\linewidth}
\centering
\begin{tikzpicture}[scale=0.75, transform shape]
\begin{axis}
[
    ybar,
    ymin=0,
    ymax=550,        
    ytick={0,100,200,300,400},
    xtick={0.2,0.4,0.6,0.7,0.8,1},
        grid=major,
         xlabel={SoC}, ylabel = {Number of vehicles},
     legend entries={Initial distribution of the SoC,
     Final distribution of the SoC}
]
\addplot +[
    hist={
     %   density,
        bins=20,
        data min=0.2,
        data max=1,
    }] table [y index=0] {initial_final_soc_case_1.txt};
    \addplot +[
    hist={
    %    density,
        bins=20,
        data min=0.2,
        data max=1
    }] table [y index=1] {initial_final_soc_case_1.txt};
\end{axis}
\end{tikzpicture}
\captionsetup{font=scriptsize}
\caption{Initial and final distribution of the SoC}
\label{final_initial_distrib}
% Attention, ici l'histogramme est construit d'une façon différente de la mienne, j'utilise dans l'article isgt les index de la méthode volume finie, ici l'approximation est différente.
\end{minipage}
\begin{minipage}[t]{0.32\linewidth}
\centering
\begin{tikzpicture}[scale=0.75, transform shape]
\begin{axis}[
    ybar,
    ymin=0,
    ymax=80,        
    ytick={0,10,20,30,40,50,60,70,80},
    xtick={1,2,3,4,5},
        grid=major,
         xlabel={Time (h)}, ylabel = {Number of vehicles}
]
\addplot +[
    hist={
     %   density,
        bins=40,
        data min=0,
        data max=5,
    }] table [y index=0] {first_time_soc_07.txt};

\end{axis}
\end{tikzpicture}
\captionsetup{font=scriptsize}
\caption{Distribution of  times PEVs SoC to attain 0.7}
\label{first_time_soc_07}
\end{minipage}
\begin{minipage}[t]{0.32\linewidth}
\centering
\begin{tikzpicture}[scale=0.75, transform shape]
      \begin{axis}[
       xlabel={Time(h)}, ylabel = {SoC},
    ymin = 0,
    ymax = 1,
    xmin = 0,
    xmax = 5,
    ytick ={0.2,0.4,0.6,0.7,0.8,1},
    grid=major,
    cycle list name=color,
        ]
        
         \addplot  [red, mark = none]    table[x index= 0, y index = 402] {soc_hist_case_1.txt};
         \addplot  [green, mark = none]    table[x index= 0, y index = 403] {soc_hist_case_1.txt};
         \addplot  [ cyan,mark = none]    table[x index= 0, y index = 404] {soc_hist_case_1.txt};
         \addplot  [magenta, mark = none]    table[x index= 0, y index = 405] {soc_hist_case_1.txt};
         \addplot  [yellow, mark = none]  table[x index= 0, y index = 406] {soc_hist_case_1.txt};
         \addplot  [ gray,mark = none]  table[x index= 0, y index = 407] {soc_hist_case_1.txt};
         \addplot  [brown, mark = none]  table[x index= 0, y index = 408] {soc_hist_case_1.txt};
         \addplot  [lime, mark = none]  table[x index= 0, y index = 409] {soc_hist_case_1.txt};
         \addplot  [olive, mark = none]  table[x index= 0, y index = 310] {soc_hist_case_1.txt};
         \addplot  [orange, mark = none]  table[x index= 0, y index = 311] {soc_hist_case_1.txt};
         \addplot  [ pink,mark = none]  table[x index= 0, y index =312] {soc_hist_case_1.txt};
         \addplot  [purple, mark = none]  table[x index= 0, y index =313] {soc_hist_case_1.txt};
         \addplot  [teal, mark = none]  table[x index= 0, y index = 314] {soc_hist_case_1.txt};
         \addplot  [ mark = none]  table[x index= 0, y index = 15] {soc_hist_case_1.txt};
         \addplot  [ violet,mark = none]  table[x index= 0, y index = 316] {soc_hist_case_1.txt};
         \addplot  [blue, mark = none]  table[x index= 0, y index =317] {soc_hist_case_1.txt};
         \addplot  [ red, mark = none]  table[x index= 0, y index = 318] {soc_hist_case_1.txt};
         \addplot  [ green,mark = none]  table[x index= 0, y index = 319] {soc_hist_case_1.txt};
         \addplot  [ cyan,mark = none]  table[x index= 0, y index = 220] {soc_hist_case_1.txt};
         \addplot  [ magenta, mark = none]  table[x index= 0, y index = 221] {soc_hist_case_1.txt};
         
        \addplot  [yellow, mark = none]  table[x index= 0, y index = 222] {soc_hist_case_1.txt};
         \addplot  [ gray,mark = none]  table[x index= 0, y index = 223] {soc_hist_case_1.txt};
         \addplot  [brown, mark = none]  table[x index= 0, y index = 224] {soc_hist_case_1.txt};
         \addplot  [lime, mark = none]  table[x index= 0, y index = 225] {soc_hist_case_1.txt};
         \addplot  [olive, mark = none]  table[x index= 0, y index = 226] {soc_hist_case_1.txt};
         \addplot  [orange, mark = none]  table[x index= 0, y index = 227] {soc_hist_case_1.txt};
         \addplot  [ pink,mark = none]  table[x index= 0, y index = 228] {soc_hist_case_1.txt};
         \addplot  [purple, mark = none]  table[x index= 0, y index = 229] {soc_hist_case_1.txt};
         \addplot  [teal, mark = none]  table[x index= 0, y index = 130] {soc_hist_case_1.txt};
         \addplot  [ mark = none]  table[x index= 0, y index = 150] {soc_hist_case_1.txt};
         \addplot  [ violet,mark = none]  table[x index= 0, y index = 131] {soc_hist_case_1.txt};
         \addplot  [blue, mark = none]  table[x index= 0, y index = 132] {soc_hist_case_1.txt};
         \addplot  [ red, mark = none]  table[x index= 0, y index = 133] {soc_hist_case_1.txt};
         \addplot  [ green,mark = none]  table[x index= 0, y index = 134] {soc_hist_case_1.txt};
         \addplot  [ cyan,mark = none]  table[x index= 0, y index = 135] {soc_hist_case_1.txt};
         \addplot  [olive, mark = none]  table[x index= 0, y index = 136] {soc_hist_case_1.txt};
         \addplot  [orange, mark = none]  table[x index= 0, y index = 137] {soc_hist_case_1.txt};
         \addplot  [ pink,mark = none]  table[x index= 0, y index = 138] {soc_hist_case_1.txt};
         \addplot  [purple, mark = none]  table[x index= 0, y index = 139]{soc_hist_case_1.txt};
         \addplot  [ cyan,mark = none]  table[x index= 0, y index = 140] {soc_hist_case_1.txt};
         \addplot  [ magenta, mark = none]  table[x index= 0, y index = 141] {soc_hist_case_1.txt};
        \addplot  [yellow, mark = none]  table[x index= 0, y index = 142] {soc_hist_case_1.txt};
         \addplot  [ gray,mark = none]  table[x index= 0, y index = 143] {soc_hist_case_1.txt};
         \addplot  [brown, mark = none]  table[x index= 0, y index =144] {soc_hist_case_1.txt};
         \addplot  [lime, mark = none]  table[x index= 0, y index = 45]{soc_hist_case_1.txt};
        \addplot  [ violet,mark = none]  table[x index= 0, y index = 146] {soc_hist_case_1.txt};
         \addplot  [blue, mark = none]  table[x index= 0, y index =147] {soc_hist_case_1.txt};
         \addplot  [ red, mark = none]  table[x index= 0, y index = 148] {soc_hist_case_1.txt};
         \addplot  [ green,mark = none]  table[x index= 0, y index = 149] {soc_hist_case_1.txt};
         \addplot  [ cyan,mark = none]  table[x index= 0, y index = 150] {soc_hist_case_1.txt};
         \end{axis}
    \end{tikzpicture}
    \captionsetup{font=scriptsize}
    \caption{SoC evolution over the time}
    \label{evolution_soc}
    \end{minipage}
\end{figure*}

\subsection{Signal tracking problem}
\label{sig_track}

In this subsection, a fleet of vehicles belonging to one company, but not geographically gathered together, is considered. The central planner performs an ancillary services, to the transport grid,  through signal tracking. This signal is a proportion of PEV in mode 1. The central planner aims at making the proportion of PEVs of his fleet in mode 1 follows closely this signal. When the fleet doesn't provide ancillary services, it has a nominal consumption, denoted by $U_{pred}$. The correction term $U_{cor}$ is a deviation of the PEVs proportion in mode 1, to contribute to the balance between electricity production and consumption. The consumption target $U_{tar}$ is defined as the sum of 
 $U_{pred}$ and $U_{cor}$.
Contrary to the case study in the previous subsection, no aggregated constraint is taken into account, and there is no final condition on the SoC.
The distance between the final distribution of the SoC and of the mode of charging of the PEVs, and a target distribution $\tilde{\mu}$, defined by the central planner, is penalized. 
It is assumed that the central planner knows in advance the target consumption curve $U_{tar}$.
The central planner minimizes the following cost:
\begin{equation}
\label{def_Jh_case_2}
\begin{array}{ll}
 J^h(\tilde{m},\tilde{E}):= & 
 
  \sum_{k=0}^{N_T-1}\sum_{\ell=0}^{N_h-1}
 \Theta_{0,1}(\tilde{m}^{k,\ell}_0,E_{0,1}^{k,\ell})h\Delta t\\
 &+
 \sum_{k=0}^{N_T-1}\sum_{\ell=0}^{N_h-1} \Theta_{1,0}(\tilde{m}^{k,\ell}_1,E_{1,0}^{k,\ell})h\Delta t
 \\
 &+
 \sum_{k=0}^{N_T-1}
 \lambda^k\left(\sum_{y = \ell}^{N_h-1}h\tilde{m}^{k,\ell}_1-U_{tar}^{k}\right)^2
 \Delta t
 
 \\
 &
 +
  \sum_{\ell = 0}^{N_h-1}\sum_{i\in\K}\beta^\ell_i \left(\tilde{\mu}^{y}_i-\tilde{m}^{N_T,y}_i\right)^2 h,
\end{array}
\end{equation}
where $\beta^\ell_i\geq 0$ for any $i$ and $\ell$, $\lambda^k>0$ for any $k$ and $a^+:=\max(0,a)$.
The functions $\Theta_{0,1}$ and $\Theta_{1,0}$ are defined by \eqref{def_theta01}. As explained in Section \ref{case 0}, this function aims at penalizing large transfers and frequent transfers from one mode to another one, in order to avoid the synchronization effects and the aging of the batteries.
The central planner assumes the system service through this aggregative cost: $\textstyle \sum_{k=0}^{N_T-1}
 \lambda^k\left(\sum_{y = \ell}^{N_h-1}h\tilde{m}^{k,\ell}_1-U_{tar}^{k}\right)^2\Delta t$. This cost incites the vehicles to follow the curve $U^k_{tar}$. 
 Finally, the final cost is $\textstyle \sum_{\ell = 0}^{N_h-1}\sum_{i\in\K}\beta^\ell_i \left(\tilde{\mu}^{y}_i-\tilde{m}^{N_T,y}_i\right)^2 h$. 
The optimization problem is:
\begin{equation}
\label{disc:pb:2}
\begin{array}{l}
 \underset{\tilde{m}\in\mathbb{R}^{N_m},\tilde{E}\in\mathbb{R}^{N_E}}{\inf}J^h(\tilde{m},\tilde{E}) \\
 \mbox{subject to }\, \eqref{linear:cons}-\eqref{dis_smin},
\end{array}
\end{equation}
where $J^h$ is defined by \eqref{def_Jh_case_2}. 
The final target distribution $\tilde{\mu}_0$ in this example is the uniform distribution between 0.6 and 0.8. The distribution $\tilde{\mu}_1$ is null, which means that no PEVs are desired in mode 1 at the end of the period. The nominal proportion $U_{pred}$ is shown in blue in Figure \ref{conso_evo} and the target consumption $U_{tar}$  in red. The parameters for the final cost $\beta_i^\ell$ is taken to be 50 for any $i\in \{0,1\}$ and $\ell\in \{0,\ldots,N_h-1\}$. For any $k\in\{0,\ldots, N_T-1\}$, the coefficient $\lambda^k$ for the signal tracking cost is also taken to be 50. The parameters for the switching costs are $\theta_{1,0}=\theta_{0,1} = 0.004$ and $\tilde{\theta}_{1,0}=\tilde{\theta}_{0,1} =  2$. A population of 1000 PEVs is considered. The Chambolle-Pock Algorithm is performed with parameters $\theta = \tau = \gamma = 0.5$. After 15000 iterations, an approximation of the optimal control is obtained. This control is applied according to the procedure described in Section \ref{App_alpha_method}. The proportion of the finite  population of PEVs in mode 1, over the period, is displayed in green in Figure \ref{conso_evo}, while  the proportion for a continuum of PEV is in black. We first observe that the two curves are very close, which justifies the mean field assumption. One can notice that the proportion of PEVs charging in mode 1 at each time step is closer to the target than the nominal proportion (in blue). We deduce that the optimal control allows to follow closely the target. The initial and final distributions of the population of PEVs are plotted in Figure \ref{final_initial_distrib_2}. One can observe that the final distribution is similar to the uniform distribution. Also, there is no PEV with a SoC lower than $0.6$ at the end of the period. The penalization on the final constraint enables to reach a final distribution close to the target distribution.
The number of transfers between different modes is 2278. Only 3 PEVs have never been charged at all, while 855 vehicles have 2 transfers, and 142 vehicles have 4 transfers from one mode to another. Asimple heuristic can be used in a post processing step to take over these 3 PEVs.
The average number of transfers per vehicle is 2.28. It is larger than the one in the first case study (1.85), but still satisfying. This difference can be explained by the choice of the parameters: the transfers are less penalized in the second case than in the first one, and the penalization of the deviation from the target, is significant compared to the penalization of transfers.
\section{Conclusion}
In this paper, we propose to control a large fleet of PEVs in the mean field limit framework. We formulate the problem as an optimal control of PDEs that we discretize. We prove that the discretized problem has a solution and we give sufficent condition to solve it numerically with the Chambolle-Pock Algorithm. We detail the implementation of the mean field control to a finite population of PEVs. We illustrate our method in two case studies, where we show that the optimal mean field control gives satisfying results on finite population of PEVs.

\begin{figure*}[t]
\centering
\begin{minipage}[t]{0.32\linewidth}
\centering
\begin{tikzpicture}[scale=0.75, transform shape]
      \begin{axis}[
       xlabel={Time(h)}, ylabel = {Proportion of vehicles},
       legend entries=
       {$U_{pred}$
       ,$U_{tar}$
       ,$U_{cont}$
       ,$U_{emp}$
       },
    ymin = 0,
    ymax = 0.43,
    xmin = 0,
    xmax = 5,
    ytick ={0.1,0.2,0.3,0.4,0.5,0.6,0.7,0.8},
    grid=major,
    cycle list name=color,
        ]
        
        % add a plot from table; you select the columns by using the actual name in
        % the .csv file (on top)
       \addplot [thick,blue, mark = none] table[x index= 0, y index = 1] {consumption_data_case_2.txt};
         \addplot  [thick,red, mark = none]    table[x index= 0, y index = 2] {consumption_data_case_2.txt};
        \addplot  [thick, black,mark = none]    table[x index= 0, y index = 4] {consumption_data_case_2.txt};
         \addplot  [thick,green, mark = none]    table[x index= 0, y index = 3] {consumption_data_case_2.txt};

         \end{axis}
    \end{tikzpicture}
    \captionsetup{font=scriptsize}
    \caption{Evolution of proportion of vehicles in mode 1}
    \label{conso_evo}
    \end{minipage}
    \begin{minipage}[t]{0.32\linewidth}
\centering
\begin{tikzpicture}[scale=0.75, transform shape]
\begin{axis}[
    ybar,
    ymin=0,
    ymax=950,        
    ytick={0,250,500,750,1000},
    xtick={0.1,0.2,0.3,0.4,0.5,0.6,0.7,0.8,0.9,1},
        grid=major,
         xlabel={SoC}, ylabel = {Number of vehicles},
     legend entries={Initial distribution of the SoC,
     Final distribution of the SoC}
]
\addplot +[
    hist={
     %   density,
        bins=20,
        data min=0.2,
        data max=1,
    }] table [y index=0] {initial_final_soc_case_2.txt};
    \addplot +[
    hist={
    %    density,
        bins=20,
        data min=0.2,
        data max=1
    }] table [y index=1] {initial_final_soc_case_2.txt};
\end{axis}
\end{tikzpicture}
    \captionsetup{font=scriptsize}
\caption{Initial and final distribution of the SoC}
\label{final_initial_distrib_2}
    \end{minipage}
\end{figure*}

\bibliographystyle{IEEEtran}
\bibliography{biblio}

\appendix
\section{Appendix}\label{section_appendix}

\subsection{Proof of Proposition \ref{prop_num_schemm}}
By induction, one can show that $\tilde{m}\geq 0$. 
Suppose for a $k\in \{0,\ldots,N_T-1\}$, we have for any $\ell \in \{0,\ldots,N_h-1\}$  $\tilde{m}_i^{k,\ell}\geq 0$. Then \eqref{cons_E} %and \eqref{borne_E} 
gives $\tilde{m}_i^{k+\frac{1}{2},\ell}\geq 0$. Using \eqref{pos_evo} and \eqref{neg_evo}, one can show that:
\begin{equation*}
\tilde{m}_i^{k+1,\ell}\geq \tilde{m}_i^{k+\frac{1}{2},\ell}\left(1-\frac{\Delta t}{h}(b_i(x_{\ell +\frac{1}{2}})^++b_i(x_{\ell -\frac{1}{2}})^-)\right)\geq 0.
\end{equation*}

Similarly, the conservation of the mass can be proved by induction.
For any $k\in \{1,\ldots,N_T\}$ and $i\in \mathcal{K}$, assuming without loss of generality that $b_i\geq 0$, we have from \eqref{pos_evo}:
\begin{equation*}
\begin{array}{l}
\sum_{\ell=0}^{N_h-1} 
\tilde{m}_i^{k+1,\ell} 
=  \frac{\Delta t}{h}\sum_{\ell=0}^{N_h-1}  b_i(x_{\ell+\frac{1}{2}})\tilde{m}_{k+\frac{1}{2}}^{\ell,i}
\\ \quad-\frac{\Delta t}{h}\sum_{\ell=0}^{N_h-1}b_i(x_{\ell-\frac{1}{2}})\tilde{m}_i^{k+\frac{1}{2},\ell-1}+  \sum_{\ell=0}^{N_h-1} \tilde{m}_i^{k+\frac{1}{2},\ell} \vspace{2mm}  \\
=\sum_{\ell=0}^{N_h-1} \tilde{m}_i^{k+\frac{1}{2},\ell}
+\frac{\Delta t}{h} b^{N_h-\frac{1}{2},i}\tilde{m}_{k+\frac{1}{2}}^{N_h-1,i}

=\sum_{\ell=0}^{N_h-1} \tilde{m}^{k+\frac{1}{2},\ell}_i.
\end{array}
\end{equation*}
where we used that  $b_i(x_{N_h-\frac{1}{2}})=b_i(1)=0$. Using \eqref{switch:alpha}, we have:
\begin{equation*}
\sum_{\ell=0}^{N_h-1} \tilde{m}^{k+\frac{1}{2},\ell}_i = 
\sum_{\ell=0}^{N_h-1} \tilde{m}^{k,\ell}_i.
\end{equation*}
\hfill{$\square$}

\subsection{Sufficient condition for Slater’s constraint qualification}\label{slater_suff_cond}

In this subsection, we  give a sufficient condition to show that there exists a point $(\tilde{m},\tilde{E})$ such that $\varphi(\tilde{m},\tilde{E})+ \psi(\tilde{m},\tilde{E})$ is finite, and $(\tilde{m},\tilde{E})$ is strictly feasible (i.e. inequalities \eqref{contraintes:etat}, \eqref{dis_smin} and \eqref{cons_E} are strict).
For any $i\in \mathcal{K}$, we define $\bar{d}_i$, $\underline{d}_i$ and $p_i$  by $\bar{d}_i:=\inf_{t\in [0,T]}\,\bar{D}_i(t)$, $\underline{d}_i:=\sup_{t\in [0,T]}\,\underline{D}_i(t)$ and $p_i:=\frac{\Delta t}{h}b_i(\underline{s})$.
\begin{proposition}
If, for any $i\in \mathcal{K}$, we have $\underline{d}_i=0$, there exists $i_0\in I$ such that $b_{i_0} =0$, and $\bar{d}_{i_0} = 1$ and there exist $j\in I$ and $e$ such that $b_j\geq 0$, $0<e<\bar{d}_j$ and:
\begin{equation}\label{binomial_constraint}
\sum_{k = 0}^{\underline{s}^h} {\tau_j\choose k}
p^k(1-p)^{ \tau_j  -k}< \varepsilon-e,
\end{equation}
where $\textstyle \tau_j:=
\left \lfloor \frac{N_T}{\rho_i}\right \rfloor$ and $\textstyle \rho_i:=\left\lceil \frac{1}{\bar{d}_i-e}\right\rceil\leq T$,
then the strong duality holds for problem \eqref{prob_ap:prim}.
\end{proposition}
\textit{Proof:} Given an initial condition $\bar{m}^0$, the goal is to find a control $\tilde{E}$, which strictly satisfies \eqref{cons_E}, and ensures that $\tilde{m}$ strictly satisfies \eqref{contraintes:etat}, \eqref{dis_smin}. 

The value of $\tau_j$ is defined such that, if at every $\tau_j\Delta t$, all the PEVs in mode $j$ move to mode ${i_0}$, and a proportion of $\bar{d}_j-e$ of the total population of PEV  moves to mode $j$, then all the PEV have been at least once in mode $j$ at the end of the period $T$. 

We define $\textstyle \mu_j:=\sum_{k = 0}^{\underline{s}^h} {\tau_j\choose k}
p^k(1-p)^{ \tau_j  -k}$. This quantity represents the proportion of PEVs with a SoC lower or equal to $\bar{s}^h$ after period of $\tau_j \Delta t$, when all the PEV charge in mode $j$, with a charging rate equals to  $b_j(\underline{s})$ and has a SoC equals to $0$ at the beginning of the period.
Thus, if no transfers occur between instants $r\Delta t$ and $(r+\tau_j)\Delta_j$,  one has that $\textstyle \mu_j \sum_{l = 0}^{s_h}\tilde{m}_j^{r,\ell}h$ is larger than $\textstyle \sum_{l = 0}^{s_h}\tilde{m}_j^{r+\tau_j,\ell}h$.  In other words, the quantity $\textstyle \mu_j \sum_{l = 0}^{s_h}\tilde{m}_j^{r,\ell}h$ is an upper bound of the proportion of PEVs with a SoC lower than $\underline{s}^h$ in mode $j$ after a period of $\tau_j\Delta t$. This observation comes from the definition of $\mu_j$, the fact that $b_j$ is a decreasing function, and that PEVs in mode $j$ at time $r\Delta t$ may have a SoC larger than $0$.

If \eqref{binomial_constraint} is satisfied, we consider the following procedure leading to a couple $(\tilde{m}, \tilde{E})$ which strictly satisfies constraints  \eqref{cons_E}, \eqref{contraintes:etat} and \eqref{dis_smin}. 
We consider $\rho_j$ periods of time $\tau_j \Delta t$. At the $k^{th}$ period, a proportion of $\bar{d}_j-e$ of PEVs moves to mode $j$ from the other modes. PEVs with a SoC lower than $\underline{s}^h$ have priority. Simultaneously, a proportion of $\textstyle \bar{d}_j-e-\frac{ke}{\rho_j} $ PEVs moves from mode $j$ to mode $i_0$. The proportions of PEVs to transfers from one mode to another, are such that the values of $\tilde{E}$ to perform the transfers, can be taken such that \eqref{cons_E} is strictly satisfied. In addition, the proportion of PEVs in mode $j$ over the $k^{th}$ period is equal to $\textstyle\bar{d}_j -e + \frac{ke}{\rho_j+2}$, making the constraint  \eqref{contraintes:etat} strictly satisfied. Between two periods of times, if  $\tilde{m}$ is not null, then $\tilde{E}$ can't be also null, to strictly satisfy \eqref{cons_E}. Thus, we consider extremely low transfers where $\tilde{m}$ is not equal to 0, such that the total sum of these transfers during the overall period does not exceed $\textstyle \frac{e}{\rho_j+2}$.

Using the definition of $\mu_j$ and its properties, at the end of the $k^{th}$ period we have $\textstyle \mu_j \sum_{l = 0}^{s_h}\tilde{m}_j^{k \tau_j,\ell}\geq  \sum_{l = 0}^{s_h}\tilde{m}_j^{(k+1)\tau_j,\ell}h$. Since, over the $k^{th}$ period, there is a proportion of $\textstyle \bar{d}_j-e+\frac{ke}{\rho_j+2}$ PEVs in mode $j$, we have $\textstyle \mu_j\left(\bar{d}_j-e+\frac{ke}{\rho_j+2}\right)\geq  \sum_{l = 0}^{s_h}\tilde{m}_j^{(k+1)\tau_j,\ell}h$.
Due to the constraint \eqref{cons_E}, not all the PEVs with a SoC lower than $\underline{s}^h$ can move at least once to mode $j$ to charge. To strictly satisfy this constraint, assume that at most a proportion of $\textstyle \frac{e}{\rho_j +2}$ PEVs, with a SoC lower than $\underline{s}^h$, does not visit at all the mode $j$. 
For the last period, only a proportion of PEVs of $\textstyle \left(1-(\bar{d_j}-e)\rho_j)\right)^+$ moves to mode $j$. At the end of the last period, the proportion of PEVs with a SoC lower than $\underline{s}^h$ is bounded above by $\mu_j+e$. Using the inequality \eqref{binomial_constraint}, the constraint \eqref{dis_smin} is strictly satisfied. We deduce that $(\tilde{m},\tilde{E})$ satisfies the Slater’s constraint qualification condition. In addition, $(\tilde{m},\tilde{E})$ is in the domain of $J^h$. 
\hfill{$\square$}

\end{document}